\documentclass[11pt]{amsart}
\usepackage{amssymb,amsmath}
\usepackage{epsf,latexsym,graphicx}
\usepackage{psfrag}
\psfragscanon
\theoremstyle{plain}
\newtheorem{thm}{Theorem}[section]
\newtheorem{prop}[thm]{Proposition}
\newtheorem{cor}[thm]{Corollary}
\newtheorem{lemma}[thm]{Lemma}

\theoremstyle{definition}

\numberwithin{equation}{section}

\newcommand{\oo}{\omega}
\newcommand{\om}{\boldsymbol{\omega}} 

\title{Differential Operators and Weighted Isobaric Polynomials}
\author{Trueman MacHenry and Geanina Tudose}
\begin{document}
\thispagestyle{empty}
\begin{abstract}
We characterize those sequences of weighted isobaric polynomials [5] which
belong to the kernel of the linear operator $D_{11} - \sum_{j=1}^k a_j t_j
D_{2j} - mD_2$,  and we characterize those linear operators of
this form in terms of the coefficients $a_j$ which have a non-zero
kernel.
\end{abstract}

\address{Department of Mathematics and Statistics\\
York University\\
North York, Ont., M3J 1P3\\
CANADA and School of Mathematics, University of Minnesota, Minneapolis, MN, 55455, USA} 
\email{machenry@mathstat.yorku.ca}
\email{tudose@math.umn.edu}
\maketitle

\footnotetext[3]{ G. Tudose supported in part by NSERC. }


\section{Introduction}
  In [4]  the following linear operator was introduced:
$$
\mathfrak{T}_m = D_{11} -\sum_j  t_j D_{2j} - mD_2,
$$
\noindent 
where   $m\in \mathbb{Z}$.                            
We are interested in these operators as linear operators on a special  
ring of polynomials, discussed in [5], namely, the ring of
{\it isobaric polynomials} $\tilde{\Lambda}$, a ring isomorphic to the
ring   of symmetric functions $\Lambda$.  The polynomials in
$\tilde{\Lambda}$, or more precisely $\tilde{\Lambda}_k$,
are over indeterminants $t_1, \ldots, t_k$  and the isomorphism just
mentioned is given by identifying $t_j$  with 
the signed elementary symmetric polynomial  $(-1)^{j+1}e_j$.  This determines 
an involutory mapping whose elements in  $\tilde{\Lambda}$  are called
isobaric reflects.

  Isobaric polynomials can be defined independently of  $\Lambda$  as
follows: for each $n$,  let $P_{k,n}: = \sum_\alpha A_\alpha
t^\alpha$    where $\alpha = (\alpha_1,\ldots , \alpha_ k)$,
$t^\alpha= t_1^{\alpha_1} \ldots t_k^{\alpha_k}$, and   $
(1^{\alpha_1}, \ldots, k^{\alpha_k})$ is a partition of  $n$, $n$  
remaining constant for all monomials in the polynomial. (We call this a 
polynomial of degree $n$, even though the coefficient of the term of 
(maximal) degree $n$ may be zero. In general the   $A_\alpha$   can be taken from any 
commutative ring.  In what follows we restrict ourselves  to the ring of 
integers.  
  It is our purpose in this paper to discuss the portion of the 
kernel of  $\mathfrak{T}_m$    lying in $\tilde{\Lambda}$, or rather, more precisely, we are 
interested in certain sequences of polynomials in $\tilde{\Lambda}$  all of whose
entries lie in the kernel of $\mathfrak{T}_m$.  These are the sequences of 
polynomials determined by a certain weighting operation whose entries lie 
in the module of Weighted Isobaric Polynomials (WIP-module) in
$\tilde{\Lambda}$ as  defined in [5].
   
        In [4] it was shown that two such sequences are the sequence of 
Generalized Fibonacci Polynomials $\{F_n\}_n$, for $m = 2$, and the sequence of 
Generalized Lucas Polynomials $\{G_n\}_n$, for $m = 1$, which are, respectively, 
reflects of the complete symmetric polynomials (CSP) and of the power 
symmetric polynomials (PSP)(see also [1],[2],[3]. In this paper  we shall 
answer the following 
questions.  
\begin{enumerate}
\item    Are there any other WIP sequences of solutions for these two values of   
$m$ ?

\item Are there any other WIP sequences of solutions for other integer values 
of $m$ ?

\item    Are there any reasonable generalizations of these partial 
differential equations for which there are WIP sequences of solutions ?
\end{enumerate}

\subsection{ The Differential Lattice}

        With each monomial $t_1^{\alpha_1}\ldots t_k^{\alpha_k}$ we associate a lattice
$\mathcal{L}(t )$  as follows:   the top of the lattice is the node
$t_1^{\alpha_1}\ldots t_k^{\alpha_k}$, 
an element of depth $\sum_i \alpha_i$.  The elements of depth $\sum_i
\alpha_i - j$    are those in which some monomial $t^\beta=
t_1^{\beta_1}\ldots t_k^{\beta_k}$ 
of depth   $\sum_i (\alpha_i - j +1)$   has been 
replaced by a monomial $t_1^{\beta_1}\ldots t_i^{\beta_i-1}\ldots t_k^{\beta_k}$.  Depth $1$  consists of 
nodes of monomials   $ t_1 , \ldots , t_k$.   Two nodes are connected by an 
edge if one node is derived from the other by subtracting $1$   from the 
exponent sum.  One node is less then another node if its depth is smaller 
and if the two nodes are connected by a sequence of edges.  The name 
{\it differential} lattice is appropriate because, as is obvious, the lattice is 
formed by partial differentiation (forgetting the differentiation 
constant).

        Given $n$  and the {\it exponent vector}  $(\beta_1,\ldots,\beta_ k)$, we can recover the 
underlying monomial uniquely.  We shall use this fact when we refer to the 
monomial by giving its exponent vector, and abuse language somewhat by 
calling the exponent vector, the monomial itself.
(This differential lattice induces a lattice structure on the Young 
diagrams of the partitions of $n$, one which is different from the Young 
lattice.) This is illustrated for the monomial  $t_1t_2^2t_3^2$ by  
\vskip0.1cm
\noindent {\bf Example 1.}
\vskip0.1cm  
\begin{figure}[!htp]
\centering
\includegraphics{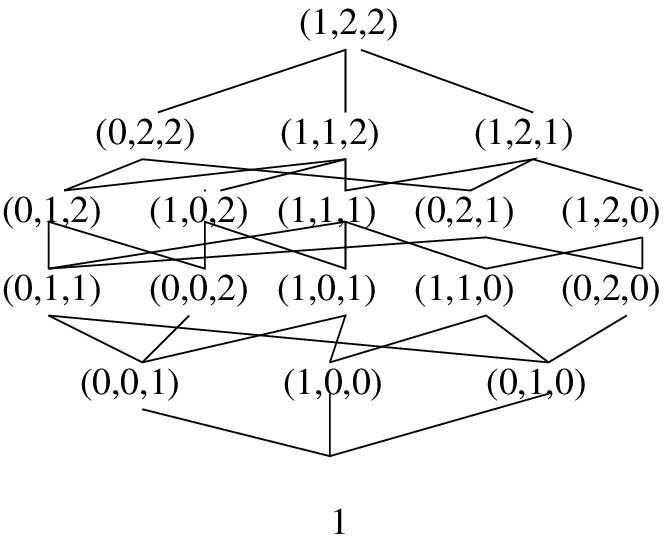}
\end{figure} 

This lattice is useful in organizing the definition of weighted monomials. 
{\it Weighting} is a method of systematically supplying coefficients to a 
sequence of polynomials indexed by the natural numbers.  Here we note that 
for any isobaric polynomial of degree  $n$   the maximal number of 
monomials is just the number of partitions of  $n$. In fact, we can regard 
an isobaric polynomial of degree   $n$   as one which is indexed by the 
shapes for the partitions of $n$.

\section{ Weighted Isobaric Polynomials}

        The concept of {\it Weighted Isobaric Polynomials} was introduced 
in [5].  We recall the definition here.  We assign a weight (here, 
integers) 
to a monomial  $t^\alpha$   by first assigning a weight $\omega_ j$ to the variable $t_j$ for $j 
= 1,2,\ldots$, where  $\om = (\omega_ 1,\omega_2, \ldots)$  is a
{\it weight vector}.    Then for any 
monomial in the lattice $\mathcal{L}(t)$ we assign the sum of weights of all nodes 
that are connected to that node by an edge.  For example, after assigning 
the weights   $\omega_j$     to the variable   $t_j$ in the previous example, the 
monomial (whose exponent vector is)   $(1,0,1)$ gets the weight
$\omega_1+ \omega_ 3$, 
while the monomial   $(0,2,0)$  gets the weight $\omega_ 2$, and, after a rather 
tedious calculation using the assignment rule, the monomial $(1,2,2)$  gets 
the weight   $6( \omega_1+2\omega_ 2+2\omega_ 3)$. Fortunately, we can avoid this calculation 
using the following theorem [5], Theorem 1).
\begin{thm}\label{thm1}
  Given a weight vector $\om = (\omega_1, \omega_2, \ldots)$  the weight  assigned 
to the monomial whose exponent vector is  $(\alpha_1, \ldots, \alpha_ k)$
is 
\begin{equation}\label{tthm1}
A_\alpha= { {\sum \alpha_j } \choose {\alpha_1 \ldots \alpha_k}} \frac{ \sum_j
\alpha_j \omega_j}{ \sum _j \alpha_j}.
\end{equation}      
\end{thm} 

        Thus any weighted isobaric polynomial is of the form
$P_{n,\om}   = \sum_\alpha  A_\alpha t^\alpha$,  
where  $\alpha$  ranges over the partitions of $n$. Each weight vector
determines a unique sequence of WIPs.  Two such sequences are
$\{F_n\}$  and $\{G_n\}$; the 
coefficients for the $F$-sequence are given by $\displaystyle{
A_\alpha =  { {\sum \alpha_j } \choose {\alpha_j}}   }$  with 
weight vector    $\om =(1,\ldots ,1, \ldots)$, i.e., $\omega_j = 1$  for all 
$j$, and for the $G$-sequence  $\displaystyle{
A_\alpha =  \frac{ (\sum_j \alpha_j -1)!}{\prod_j \alpha_j!} n   }$,
with weight vector $\om= (1,2,\ldots ,k,\ldots)$ , i.e., $\omega_j
= j$, for all $j$.  It is easily seen that 
this assignment follows from Theorem 2.1.  We call the weight vector of 
the  
$F$-sequence, the {\it unit weight vector}, and the the $G$-sequence,
{\it the natural weight vector}.    

        In [5], Theorem 2.3, it is shown that the sequences of weighted 
isobaric polynomials form a free $\mathbb{Z}$-module where addition is
defined as addition  of weight vectors, that is the sum of two
sequences of weights $\om$ and $\om'$      
is the sequence of weight $\om'' = \om  +\om'$.  It is also shown in that 
same paper that  isobaric reflects of hook Schur polynomials (i.e., the 
Schur polynomials determined by hook Young diagrams) are in the  
WIP-module. (The weight of the hook reflect determined by the hook diagram 
$(n-r,1^r)$  is  $(-1)^r(0,\ldots,0,1,1,\ldots)$ ).   The hook reflects in fact form a basis for the WIP- module.  As an application of the WIP-module 
structure we have the following isobaric version of a well-known theorem 
of symmetric functions
\begin{thm}
$ G_n = (-1)^r H_r$, where   $H_r$  is the hook reflect induced by the 
shape $(n-r,1^r)$.
\end{thm}
For, clearly the sum of the alternating sum of the $n$-hook weights is
the weight of $G_n$.
\hfill
$\Box$

The symmetric polynomial version of this is the statement that a
complete symmetric polynomial is an alternating sum of Schur hooks.

\section{The kernel of $\mathfrak{T}_m$}

We now turn our attention to the linear operator $\mathfrak{T}_m$ and find that for
certain choices of the parameter $m$, the $F$-sequence and the
$G$-sequence belong to the kernel of $\mathfrak{T}_m $.

\begin{thm}[\cite{TM3},Theorem 4]\label{old}
$ \mathfrak{T}_m(F_n)  = (D_{11} -\sum_j  t_j D_{2j} - mD_2)(F_n) = 0$
when   $m = 2$, and  $ \mathfrak{T}_m(G_n)  = (D_{11} -\sum_j  t_j
D_{2j} - mD_2)(G_n) = 0$ when $m=1$.
\end{thm} 
        
This theorem will follow from Theorem 3.2 below. Theorem
3.1 tells us that the $F$- and $G$-sequences are  
solutions to the operator equation when the parameter is $m = 1$ in the 
case of the $G$-polynomials and $m=2$ in the case of the $F$-polynomials, but it turns out that these solutions are determined by other 
more basic solutions which, while dependent on the weights of the $F$- and
$G$-sequences, are not themselves WIPs. We might refer to these
polynomials as {\it satellites}. We call them {\it strings}.

  Thus, for each
 $F_n$  and $m=2$  there is a sequence of isobaric polynomials \\
$S_1^{F_n},\ldots, S_u^{F_n}$,  where $F_n =\sum_{j=1}^u S_j^{F_n}$
and  $\mathfrak{T}_2( S_j^{F_n}) = 0$.       
Similarly, for each  $G_n$, and $m = 1$, there is a 
sequence of isobaric polynomials    
$ S_1^{G_n}, \ldots, S_u^{G_n}$,  where $G_n= \sum_{j=1}^u S_j^{G_n}$
and  $\mathfrak{T}_1( S_j^{G_n}) = 0$.

        To see this,   for a given $n$ and $k > 1$, we first choose
        exponent vectors  of the following two kinds:
\begin{enumerate}
    
\item vectors of type $(0, \alpha_2, \alpha_3,\ldots , \alpha_k)$,
where    $\alpha_3,\ldots ,\alpha_k$  is a fixed 
$(k-2)$-tuple and   $\alpha_2$  is largest second element with respect
to this condition,
   
\item vectors of type  $(1, \alpha_2,\alpha_ 3,\ldots , \alpha_k)$
  where  $\alpha_3,\ldots ,\alpha_k$  is a fixed 
$(k-2)$-tuple and   $\alpha_2$  is largest second element with respect
to this condition.
\end{enumerate}
  
Then we select sequences of these vectors of the following form
$$
\begin{aligned}
(0,\alpha_2,\alpha_3,\ldots,\alpha_k) &\qquad  (1,\alpha_2,\alpha_3,\ldots,\alpha_ k)\\
(2,\alpha_2-1,\alpha_3,\ldots,\alpha_k)  &
                         \qquad (3,\alpha_2-1,\alpha_3,\ldots,\alpha_ k)\\
\ldots &   \qquad    \ldots\\
(2j,\alpha_2-j,\alpha_3,\ldots,\alpha_k)  & \qquad (2j+1,\alpha_2-j,\alpha_3,\ldots,\alpha_ k)\\  
\ldots & \qquad      \ldots\\  
(2\alpha_2,0,\alpha_3,\ldots,\alpha_k)  & \qquad (2\alpha_2+1,0,\alpha_3,\ldots,\alpha_ k)
\notag
\end{aligned}   
$$  
$$
                                 j = 0, 1, \ldots ,\alpha_2.
$$
Such a sequence is called {\it a string}.  The  first element in the sequence is the 
{\it string generator}. If the string starts with $0$, we call it an
{\it even 
string}, and if it starts with $1$, we call it an {\it odd string}. The 
left-hand column above is an {\it even string},  while the right-hand 
column is an {\it odd string}.  It is not 
difficult to see that for a given   $n$   all of the exponent vectors that 
arise from the partition of   $n$   occur in some even or odd string.  Thus, 
every isobaric polynomial is just 
 the sum of its strings with ``remembered''
coefficients.  In particular, for a sequence of weighted isobaric 
polynomials, each polynomial is just the weighted sum of its strings. 
Theorem 3.1 will follow from this fact. We shall say that a 
string {\it belongs} 
to a weighted isobaric polynomial if it is a weighted string in that 
polynomial. For example, the (three) strings that {\it belong} to $F_4$ , 
where  $F_4 = t_1^4+3t_1^2t_2+t_2^2+2t_1t_3+t_4$, are
(0,2,0,0), (2,1,0,0), (4,0,0,0); (1,0,1,0); and  (0,0,0,1).

It is clear that the strings reflect the truncations of the isobaric 
polynomial obtained by deleting the variables $t_j$ for the $j$'s from a 
certain $j$ onward.

\begin{thm}\label{string}\hfill

 (1)  If $S^F$  is a string belonging to $F$, then   $\mathfrak{T}_2(S^F) = 0$.

(2)   If   $S^G$  is a string belonging to $ G$, then $\mathfrak{T}_1(S^G) = 0$.
\end{thm}

This theorem will follow from
 
\begin{lemma}\label{l1}

a).  $(2j+2,  \alpha_2-j-1, \alpha_3, \ldots, \alpha_k)$ and $(2j,
  \alpha_2-j,  \alpha_3, \ldots, \alpha_k)$  are 
adjacent elements in the even string generated by $(0, \alpha_2,
  \alpha_3,\ldots, \alpha_k)$  and 
the coefficient of   $D_{11}(2j+2, \alpha_2-j-1,  \alpha_3, \ldots ,
  \alpha_k)  = -(\mathfrak{T}_ m - D_{11})(2j,  \alpha_2-j,  \alpha_
  3, \ldots,  \alpha_k)$ 
whenever the weight vector is $(1,1,\ldots,1,\ldots)$  and   $m = 2$, or the weight 
vector is $(1,2,\ldots,k,\ldots )$  and $m = 1$.  $D_{11}$  applied to the string 
generator is $0$   and  $(\mathfrak{T}_ m - D_{11})$ applied to the
  last element in the string is also   $0$.
 
b).  $(2j+3,  \alpha_2-j-1, \alpha_3, \ldots, \alpha_k)$ and $(2j+1,
  \alpha_2-j,  \alpha_3, \ldots, \alpha_k)$  are 
adjacent elements in the odd string generated by $(1, \alpha_2,
  \alpha_3,\ldots, \alpha_k)$  and 
the coefficient of   $D_{11}(2j+3, \alpha_2-j-1,  \alpha_3, \ldots ,
  \alpha_k)  = -(\mathfrak{T}_m - D_{11})(2j+1,  \alpha_2-j,  \alpha_
  3, \ldots,  \alpha_k)$ 
whenever the weight vector is $(1,1,\ldots,1,\ldots)$  and   $m = 2$, or the weight 
vector is $(1,2,\ldots,k,\ldots )$  and $m = 1$.  $D_{11}$  applied to the string 
generator is $0$   and  $(\mathfrak{T}_ m - D_{11})$ applied to the
  last element in the string is also   $0$.
\end{lemma}
 
\noindent
{\bf Proof of Lemma} 
That the elements mentioned in the lemma belong to the 
string and are adjacent is obvious.  The fact that the first and last
elements of the string are mapped to $0$  by the operators   $D_{11}$    and 
$(\mathfrak{T}_m - D_{11})$ as claimed is also obvious.  We shall prove then that the 
coefficients of the elements $D_{11}(2j+2, \alpha_2-j-1,  \alpha_3, \ldots ,
  \alpha_k)$    and $(\mathfrak{T}_ m - D_{11})(2j,  \alpha_2-j,  \alpha_
  3, \ldots,  \alpha_k)$ ultiplying
are negatives of one another.
        
By Theorem 2.1 we have that
\begin{equation}\label{1}
A_{s_{2j+2}}=
\frac{(\sum_{i=2}^k \alpha_i+j)!}{
  (2j+2)!(\alpha_2-j-1)! \prod_{i=3}^k\alpha_i!}
[(2j+2)\omega_1+(\alpha_2-j-1)\omega_2+\sum_{i=3}^k\alpha_i\omega_i]
\end{equation}
\begin{equation}\label{2}
A_{s_{2j}}=
\frac{(\sum_{i=2}^k \alpha_i+j-1)!}{
  (2j)!(\alpha_2-j)!\prod_{i=3}^k\alpha_i!}
[(2j)\omega_1+(\alpha_2-j)\omega_2+\sum_{i=3}^k\alpha_i\omega_i],
\end{equation}
\noindent
where $s_{2j+2}=(2j+2,\alpha_2-j-1,\ldots, \alpha_k) $ and $s_{2j}=(2j,\alpha_2-j,\ldots, \alpha_k)$.
The coefficient due to $D_{11}$  applied to $s_{2j+2}$  is
\begin{equation}\label{3}
  (2j+2)(2j+1),
\end{equation}
\noindent
and the coefficient due to $\mathfrak{T}_m - D_{11}$ applied to
$s_{2j}$ is
\begin{equation}\label{4} 
          (\alpha_2-j)( \sum_{i=2}^k\alpha_i+j+m-1).
\end{equation}
Multiplying equation (3.1) by (3.3) and equation (3.2)
by (3.4) and using the values given 
by the hypothesis of the lemma for $m$ and for the weight vector and 
comparing gives the result.   

It is useful to record the last steps in the computation beginning
just before the hypotheses on $m$ and the weight vectors are applied. We 
have
this expression

\begin{equation}\label{ex1}
(\sum_{i=2}^k\alpha_i+j)((2j+2)\omega_1+(\alpha_2-j-1)\omega_2+\sum_{i=3}^k\alpha_i\omega_i)
\end{equation}
$$
-(\sum_{i=1}^k\alpha_i+j+m-1)((2j)\omega_1+(\alpha_2-j)\omega_2 +\sum_{i=3}^k \alpha_i\omega_i)
$$
Letting $m=1$  gives $2\omega_1-\omega_2=0$, after applying the
hypothesis on the weights, which gives the result required no matter
what the weights $\omega_j$ are for $j>2$. Thus we have proved more in
this case, that is, we have infinitely many WIP sequences as solutions.

Letting $m=2$ gives the expression
\begin{equation}\label{ex2}
(\sum_{i=1}^k\alpha_i+j)(2\omega_1-\omega_2)- ((2j)\omega_1+(\alpha_2-j)\omega_2+\sum_{i=3}^k\alpha_i\omega_i)
\end{equation}
\noindent
but now we need our weight hypothesis on all of the weights to achieve
the cancellation, thus this expression is $0$ if we assume that
$\omega_j=\omega_1$, for all $j$.

The proof of part b). is similar to that of part a). and will be omitted.
\hfill
$\Box$
\vskip10pt
But then Theorem 3.2 now follows from the
lemmas. Theorem 3.1 follows from  Theorem 3.2 by linearity.
\vskip10pt
        
It is an interesting consequence of the proof that the lattices of 
the string elements intersect for the first time exactly at the nodes 
determined by   $D_{11}$ operating on the string.  We give an example.  
\vskip0.1cm
\noindent {\bf Example 2.} {\it Consider the string generated by (0,2,1), 
n = 7, k = 3.  The lattices are given by}
\begin{figure}[!htp]
\centering
\includegraphics{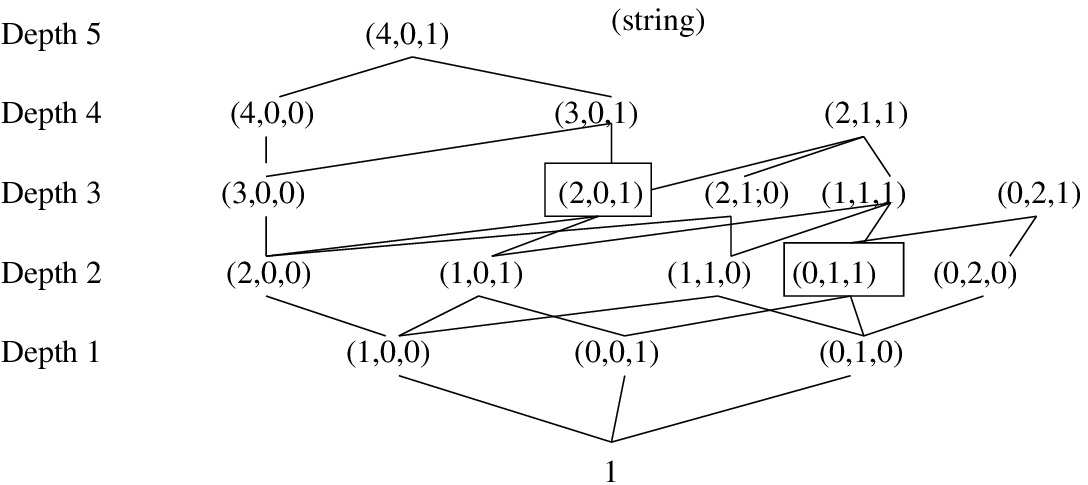}
\end{figure} 
\vskip0.1cm
In this case the intersection nodes are   $(2,0,1)$ and $(0,1,1)$. The
string consists of the three nodes $(4,0,1)$, $(2,1,1)$ and the string 
generator is  $(0,2,1)$.

   It is also the case that the intersection nodes again form a 
string. This time for polynomials of degree   $n-2$.  However, these
strings do not inherit the weighting of the string of degree $n$.
        
\section{Other solutions}

\indent We stress here that what we mean by a solution is the entire 
sequence of WIPs determined by a particular weight vector $\om$; 
calling such a sequence  $P_{\om}$,  $P_{\om} = \lbrace 
P_{n,\om}\rbrace$,we have as solutions the polynomials generated by the 
strings of WIP-solutions.  We claim that the  WIP sequences of solutions 
of 
the PDE $\ \mathfrak{T}_ 1 = 0$ (that is, for $m = 1$)  are exactly those
solutions generated by linearity from  
the strings in which $2\omega_1 = \omega_2$ with the $\omega_j$  
arbitrary for $j>2$,  but fixed 
throughout the string;  $G$-polynomials, for example. 

        When $m = 2$, the solutions of $\mathfrak{T}_2= 0$ consist just of the scalar 
multiples of the   $F$-polynomials. The kernel of the operator    
operating on the WIP-module is  $0$    when $m \neq   1,2 $.  

The $G$-polynomials are one of many WIP-solutions of the operator
equation when $m=1$; but when $m=2$, the only WIP solutions are just
the scalar multiples of the $F$-polynomials.
  We prove these assertions now.

\begin{prop}\label{p}
  Let $\om$    be a weight vector and  $\{P_{n,\om}\}$   be a
  sequence of solutions  of    $P_{n,\om} \in $WIP-module, then
  either

 (1) $m = 1$ and $2\omega_1 = \omega_2$ or 

 (2)  $m = 2$ and  $\omega_ 1 = \omega_2$.
\end{prop}

\noindent
{\bf Proof}
It is only necessary to look at the second and third terms of the 
sequence; namely, at                                                            $$                                        
P_{2,\om}  = \omega_1t_1^2 + \omega_ 2t_2 , 
$$
$$                       
P_{3,\om}  = \omega_1 t_1^3 +(\omega_1+ \omega_2)t_1t_2+ \omega_3t_3 .  
$$
Requiring that   $P_{2,\om}$    satisfies the operator equation
implies that  $m \omega_2= 2 \omega_1$ ; requiring that
$P_{3,\om}$   satisfies yields  $ m (\omega_1+\omega_2) = (5
\omega_1-\omega_2)$.
  
\noindent
Setting the two values equal and solving the resultant quadratic in
$\mathbb{Z}$ gives the 
two possibilities    $ 2\omega_1 = \omega_2,\  \omega_1 =   \omega_2$
or $\om = 0$.  Solving for $m$  in each case 
gives $m = 1$ and $m = 2$ or the trivial case for any $m$,
respectively.  And we know that the first  two cases are 
realized with the $G$-polynomial sequence and the $F$-polynomial sequence 
respectively.      
\hfill
$\Box$
\vskip10pt

We summarize our discussion of solutions of the operator so far in the
following Corollary of Proposition 4.1.

\begin{cor}\label{p'}
In case (1) of Proposition 4.1, the condition $2\omega_1 =  
\omega_2$  characterizes the solutions of the operator equation $\mathfrak{T}_1=0$.  In 
particular,  the first two components of the weight vector completely 
determine the kernel.
                                
In case (2) of Proposition 4.1, if we ask 
for solutions generated by string solutions then  $\omega_1 =\omega_2
=a$      implies that $\omega_j =a$ for all   $j$. That is, all such
WIP-solutions are of the form    $a  F_n$, $n\in \mathbb{N}$, and all solutions are exactly those generated in the WIP-module by $F$-strings.
\end{cor}
        
\noindent
{\bf Proof}
The proof consists of looking at the proof of Theorem 3.2 more 
carefully and noting that in light of Proposition 4.1 (1),  the 
cancellation 
occurs independently of the choice of  $\omega_j$   for $j > 2$, while
in the case of Proposition 4.2 (2), the proof arrives at the
equation  
 $\sum_j a \alpha_j = \sum_j \alpha_j \omega_j$  with   
$\omega_1 =  \omega_2 =a$, which must hold for all exponent vectors
$\alpha$   and for a fixed 
weight vector $\omega$;  thus    $\omega_j = a$  for all $j$.    
\hfill
$\Box$

\vskip 24pt

So now we come to the three questions posed in the introduction. It
turns out that we shall be able to answer these questions completely
once we have answered the third one. So our aim is to prove

\begin{thm}\label{last}
The operator $D_{11} - \sum_j a_j t_j D_{i,j} -mD_2$,  where  $a_j
 \in \mathbb{Z}$,  
has WIP-sequence solutions only when   $a_j = 1$ and   $m = 1$ or $m =
 2$, where the $a_j$ and $m$ are assumed to be arbitrary real numbers, not 
all zero
 (Though, they could be taken from any field of characteristics $0$ as
 far as the proof is concerned).
\end{thm}

The statement of this theorem makes clear what we have chosen to mean
by a generalized operator. 
 When one tries to find other second order, linear 
partial differential equations that have sufficient resemblance to the one 
at hand, the lack of  left-right symmetry among the partitions of $n$ as $n$ 
increases becomes more apparent.  This is due to the fact that $1$'s will 
appear in the decompositions of  $ n$   many times, but   $n$  itself can 
appear only once; small numbers have the advantage over big ones. This is 
reflected in the futility of trying to find new PDEs by varying the 
suffixes of the operators, $D_{i,j}$.  However, a tack that appears promising 
is to provide   $\mathfrak{T}_m$   with arbitrary (real) coefficients.
Thus, we want to ask what is the kernel  of any operator of the sort
$D_{11}-\sum_j a_j t_j D_{2,j}-m D_2$, $a_j$ arbitrary (real) scalars?
(Here we assume that the coefficient of $D_{11}$ is not 0, so we can, 
without loss of generality, assume that it is 1.) By the way, the 
resemblance of the operator equation   
$\mathfrak{T}_m = 0$   to the "Newton identity" 
satisfied by the
WIP-polynomials (see Theorem 4.1 [5]) is striking, and probably 
significant, though the
anomolous role of the  $D_2$-term is 
puzzling.

\vskip10pt
        Before we  prove the theorem we will prove some results which 
are of interest in their own right, and which will
        contribute directly to the proof of the theorem.

\begin{lemma}\label{ll2}
Let   $S$  be a string belonging to  $P_{n,\om} $, then
$D_{11}(P_{n,\om} )$   and   $(\mathfrak{T}_m-D_{11}) (P_{n,\om} ) \in $
the union of the lattices of the string elements; call this 
the lattice of the string.
\end{lemma}
 
\noindent
{\bf Proof}  This is clear from the definition of the differential
lattice.   
\hfill
$\Box$ 
\vskip10pt

\begin{thm}\label{str}
 $P_{n,\om}$   is a  solution of  $\mathfrak{T}_m = 0$   if and only if the strings belonging to  $P_{n,\om}$  are solutions.
\end{thm}

\noindent
{\bf Proof}  Clearly, since $P_{n,\om}$   is just the sum of its strings with 
remembered coefficients, we need only prove the necessity. 
         So suppose that $P_{n,\om}$  is a solution of 
 $\mathfrak{T}_m = 0$, but then the theorem  follows from
 Lemma 4.4.   
\hfill
$\Box$
\vskip10pt
        Next,  we note that the depths of elements in a string form a 
strictly monotonically increasing sequence.  It is also easy to see that  
$\mathfrak{T}_m (S)$  is also a string. From the fact that  $D_{11}$ 
  and    $\mathfrak{T}_m-D_{11}$  each have 
exactly one monomial as an image,  we  have, using Lemma 4.4,

\begin{lemma}\label{ll3}
   A string is a solution of   $\mathfrak{T}_m = 0$ if   only if the
   proof is the "domino" proof  used in the proof of Theorem 3.2.     
\hfill $\Box$

\end{lemma}

        The proof of Proposition 4.1 contains the following fact which, 
together with its proof,  also holds in the generalized operator case.

\begin{lemma}\label{ll4}
  If   $P_{2,\oo}$   satisfies the generalized operator equation, then 
$m \oo_2 = 2\oo_1$. \hfill $\Box$
\end{lemma} 
\vskip 0.5cm
\noindent
{\bf Proof of Theorem 4.3}
 Now, it is clear that if the strings of $P_{n,\om}$   satisfy
the operator equation, then so does $P_{n,\om}$.  So let us suppose, 
conversely, 
that $P_{n,\om}$     
satisfies the operator equation. We consider the ``domino''  proof used in 
the proof of Theorem 3.2.
Let  $\alpha= (\alpha_1,\alpha_2,  \alpha_3,\ldots, \alpha_k)$  be an
        arbitrary element in a string  $S$.  If  $\alpha$    
is the only element of the string, then clearly it satisfies the operator 
equation. So suppose that $\alpha$ 
    is not an element of least depth, that is, it is not  the generating 
element of the string.  In this case, there is an element   of depth one 
less than the depth of $\alpha$.  Let us suppose that 
 $D_{11} (A_{\oo} \alpha )  =   ( \mathfrak{T}_m-D_{11} )(B_{\oo} \beta)$,
        that is suppose that the ``domino'' proof applies. (Note that if
        $\alpha$   is an element of greatest depth, then
        $D_{11}(\alpha) = 0$.)  The picture looks like  this: 
\vskip2cm   
\hskip1.0cm   $(\alpha_ 1,\alpha_2,\alpha_3,...,\alpha_k) = \alpha $
\vskip0.2cm
\hskip2.2cm $\arrowvert$
\vskip0.2cm                                               
\hskip1.0cm   $(\alpha_1-1,\alpha_2,\alpha_3,...,\alpha_k) 
\hskip0.6cm (\alpha_1-2,\alpha_2+1,\alpha_3,...,\alpha_k) = \beta$ 
\vskip0.2cm  \hskip2.2cm $\arrowvert$ \hskip3.0cm $/$
\vskip0.2cm                        
\hskip1.0cm   $(\alpha_1-2,\alpha_2,\alpha_3,...,\alpha_k)$     
\vskip0.2cm

Recalling the proof of Theorem 3.2 at this point, we needed to equate the 
product of  the  coefficient of   $A_{\oo}$    of $\alpha$     and the
coefficient of   $D_{11}(\alpha)$  with the product of the coefficient
$B_{\oo}$     of  $\beta$     and   $(\mathfrak{T}_m-D_{11})(\beta )$.      
The new ingredient here is that   $(\mathfrak{T}_m-D_{11})(\beta ) =
(\alpha_2+1) (\sum_j a_j \alpha_j + m - 2a_1)$  
due to the new coefficients of $\mathfrak{T}_m-D_{11}$.

\noindent After making the calculation indicated above and allowing the 
dust to settle, this gives 
\begin{equation}\label{4.1}          
( \sum_{j=1}^k \alpha_j-1)( \sum_{j=1}^k \alpha_j\omega_j) = 
(\sum_{j=1}^k\alpha_j \omega_j + \omega_2 - 2\omega_1)(\sum_{j=1}^k a_j 
\alpha_j + m - 2a_1) \end{equation}

\noindent
as a necessary condition for the generalized operator to have a
solution. 

We shall assume throughout that  $P_n$  is not trivial, that is, that 
$\om \not=0$. Equation (4.1) can be rewritten as

\begin{equation}\label{4.2}
(\sum_{j=1}^k\alpha_j-1)(\sum_{j=1}^k\alpha_j\omega_j)-(\sum_{j=1}^k\alpha_j
\omega_j)(\sum_{j=1}^ka_j\alpha_j+m-2a_1)
\end{equation}
$$-(\omega_2-2\omega_1)(\sum_{j=1}^ka_j\alpha_j)=(\omega_2-2\omega_1)(m-2a_1)$$
\indent
 The left hand side of (4.2) depends on $ \alpha $ , which is a 
variable,  while the right hand side depends only on the choice of 
$ \om $ and the constants $a_1$ and $m$. Hence the left hand side and 
the right hand side of (4.2) are independently 0. And hence,  either $ 
\omega_2 - 2 \omega_1 = 0 $ or $ m - 2a_1 = 0 $.   
\vskip 0.1cm
\indent 
In the first case, we have then that $\omega_2 = 2 \omega_1
$ and,  by Lemma 4.8, $m = 1$. The left hand side of (4.2) becomes  
\begin{equation}\label{4.3}
(\sum_{j=1}^k\alpha_j\omega_j)
(\sum_{j=1}^k\alpha_j  -\sum_{j=1}^ka_j\alpha_j+2-2a_1)=0. 
\end{equation} 
\vskip 0.1cm
$(\sum_{j=1}^k\alpha_j\omega_j) = 0$  implies that $\om = 0$, that is,  the solution is trivial.
  
\vskip 0.1cm \noindent Thus 
$(\sum_{j=1}^k\alpha_j-\sum_{j=1}^ka_j\alpha_j+2-2a_1= 
0)$ and so we have as above that  
$\sum_{j=1}^k\alpha_j-\sum_{j=1}^ka_j\alpha_j = 0$ and $2-2a_1=
0$. And so we have that 
\begin{equation}\label{4.4}
 a_1= 1  \mbox{ and }
\sum_{j=1}^k\alpha_j=\sum_{j=1}^ka_j\alpha_j.  
\end{equation} 
From these equations we have that $a_j = 1$ for $j=1,\ldots,k$.
This is just the case of the original operator for which the 
$G$-polynomials were solutions.

In the second case, from  $m-2a_1 = 0$, we have $m = 2a_1$ and
from this and (4.1) we have 
\begin{equation}\label{4.5} 
\sum_{j=1}^k\alpha_j\omega_j(\sum_{j=1}^k\alpha_j - \sum_{j=1}^k\alpha_j 
a_j  - 1) = (2\omega_1 - \omega_2)(\sum_{j=1}^k\alpha_j-j a_j).
\end{equation}
 
Consider the monomial   $\omega_n t_n$. It follows from the 
definition of a string that $\omega_nt_n$ is a
string, or, in the case that $n=2$, is the generator of a two element 
string, so we apply Theorem 4.5.  Here   $\alpha_n =1$ and 
$\alpha_j = 0$ otherwise.  From (4.5) we arrive at 
\begin{equation}\label{4.6} \omega_na_n = (2\omega_1 - \omega_2)(a_n). 
\end{equation} 
We may assume
that $a_n \ne0$ for some $n$.  Then,  when  $n = 2$ we have that $\omega_2 
=\omega_1$  for all $n$,  so that  $\omega_n = 2\omega_1-\omega_2$ is 
constant 
for all $n$.  In particular,  when $n = 2$ we have that  $\omega_2 = 
\omega_1$, and thus $\omega_n = \omega_1$ for all $n$. So, in particular, 
when $n=2,\,  \omega_2 = \omega_1$ which, in turn, implies that $\omega_n = 
\omega_1$ for all $n$. From (4.1) it then follows that  
$\omega_1(\sum_{j=1}^k\alpha_j-1)(\sum_{j=1}^k\alpha_j) = 
(\sum_{j=1}^k\alpha_j-1)(\sum_{j=1}^ka_j\alpha_j)\omega_1$. If $ 
\om \ne0$,  that is,  if the solution is not the trivial solution, we 
have that $\sum_{j=1}^k\alpha_j = \sum_{j=1}^ka_j\alpha_j$ from which it 
follows that  $a_j = 1$ for $j = 1,\ldots,k$. Moreover,  $m = 2$;  and this is 
just the case of solutions generated by the strings of  $F$-polynomials and 
the original operator. \hfill $\Box$  
\vskip1cm
\indent$Remark:$ It is rather interesting that if $P_{\omega}$ is a 
solution 
of the 
operator equation and if  $\omega_2 = 0$, then either 
$\omega_j = 0$ for all j,  or $ m = 1$. This follows easily by applying 
Theorem 4.5, Lemma 4.7 and the assumptions to the strings generated by 
first $(0,1, \ldots ,1,0,0)$ and then the string generated by $(1,1,\ldots,1,0,0,0)$, 
with the last one being the exponent of $t_n$ in each case.
\vskip 0.1cm 
        We have then that the answer to question (3) is that the only WIP 
solutions for the generalized operator equation occur when   $a_j = 1$
 for all $j \in \mathbb{N}$, and   $m = 1$ or $2$. 
Thus the generalized operator has a zero kernel except in the case we 
started with, thus generalizing the  operator does not produce new 
solutions. Clearly, we have also answered question (2); allowing $m$  to vary 
beyond $1$ and $2$, in fact, over any field of characteristic $0$
produces no new solutions.  The answer to question  (1), we learn here, is yes and no.  If   $m =2$, 
then the answer is unique up to a scalar multiple, that is all 
WIP-solutions are scalar multiples of the $F$-sequence; but if $m=1$, then not only are 
scalar multiples of the $G$-sequence solutions, but also so is the sequence  
$P_{\om} $ anytime that
$2\oo_1 = \oo_2$, the remaining weights being arbitrary.  However, we also have 
satellite solutions, {\it the strings}, that get their life from the WIPs, but 
are not themselves WIPs.It is tempting to think that a weight vector for 
an initial string of  $P_{n,\om}$, i.e., the ``degree'' string (the string 
whose terminal element     $(n,0,\alpha_3,\ldots, \alpha_k)$ ) might be
weighted as $(\oo_1,\oo_2,0,\ldots,0,\ldots)$,   while the  $\om =
(\oo_j)$  where
 $\oo_j$  is different from $0$. The following example shows what goes
 wrong here.
\vskip0.1cm
\noindent {\bf Example 3.} {\it Consider  $ P_{n,\om}= 
\oo_1t_1^4+(2\oo_1+\oo_2)t_1^2t_2 + \oo_2t_2^2+(\oo_1+\oo_3)t_1t_3 + 
\oo_4t_4$. The strings are}: 
\begin{eqnarray*}
& &\;\>Initial \; String\\
& &{\begin{cases}&(0,2,0,0) \qquad\qquad \Big{\{} (1,0,1,0)  \qquad\qquad 
\Big{\{} (0,0,0,1)\cr
&(2,1,0,0)\cr
&(4,0,0,0)
\end{cases}}
\end{eqnarray*}
\noindent
Try weight vector $(\oo_1 ,\oo_2 ,0,0)$.  But now by Theorem 3.1, if the 
initial string is a WIP, then $(1,0,1,0)$, that is the monomial
        $t_1t_3$,  has 
coefficient   $\oo_1+\oo_3$, while $t_1t_3$  has weight  $\oo_1$  in the new 
weighting---recall, we have to assign a weight to each of the monomials 
induced by a partition of $n$, thus $t_1t_3$  would appear in the
        initial string if   $\oo_1 \neq  0$.  This contradiction would appear more generally.  We omit the proof.

        It is also interesting to note the rather special companionable 
role that the $F$-sequences and $G$-sequences play among the isobaric 
polynomials, especially among the WIPs.  In addition to the properties 
shown in this paper, we have, for example, that the $G$'s are related to the 
$F$'s by partial differentiation as follows:
 $\partial/ \partial t_j (G_n) = F_{n-j}$.  In general,  $\partial
 /\partial  t_j( P_{\om})$ is not a WIP, in fact, there 
is good reason to believe that this is the only case. We pursue this 
observation in a later paper.

\end{document}